# Light pillars over rippled water: why their sides appear parallel

Philip W. Kuchel

*School of Life and Environmental Sciences, University of Sydney, NSW 2006, Australia*

**Abstract**

Streetlights shining across rippled water often produce tall, narrow reflections with strikingly parallel sides. These 'light pillars' appear almost architectural, yet the water surface is neither vertical nor smooth. We develop a geometric-optics model that explains the phenomenon using the specular reflection rule, projection geometry, and the physics and statistics of surface slopes. A rippled water surface can be viewed as an ensemble of small facets (tangent planes on waves) acting as tiny mirrors. As one looks farther across the water (coordinate $x$), the physical width of the region whose facets reflect the light into the eye (or camera/pinhole) increases, but the pinhole projection onto the image plane compresses this widening by a factor proportional to $1/x$. The two effects cancel, producing a reflection whose image width remains constant.

This paper appears to be the first to formally document this. (Deviations from the strict cancellation at both ends of the image are also explained.)

The analysis clarifies how earlier qualitative treatments failed through a lack of identifying the role of projection geometry. The model also qualitatively explains the brightness variation along the pillar, and why extended objects such as buildings do not produce elongated reflections. The treatment is intended for scholars of introductory optics and for those interested in the physics and mathematics of everyday visual phenomena. We conclude by outlining how the same framework extends naturally to more complex viewing situations, such as the setting sun viewed from a cliff top over the sea.

## 1. Introduction

Figure 1 is one of many images on the Internet of Bratislava on the Danube at dusk; it shows the city apparently supported on tall straight-sided pillars of light projecting down from the river's edge. These 'light pillars' appear tall, narrow, and remarkably parallel-sided, as though the water were acting as a long vertical mirror. Yet the water surface is a horizontally extended, wind-rippled surface whose local slopes vary continuously in space and time. Their geometry is striking – parallel sided, almost architectural – and it immediately invites the question: Why should a natural optical phenomenon produce columns with such parallelism?

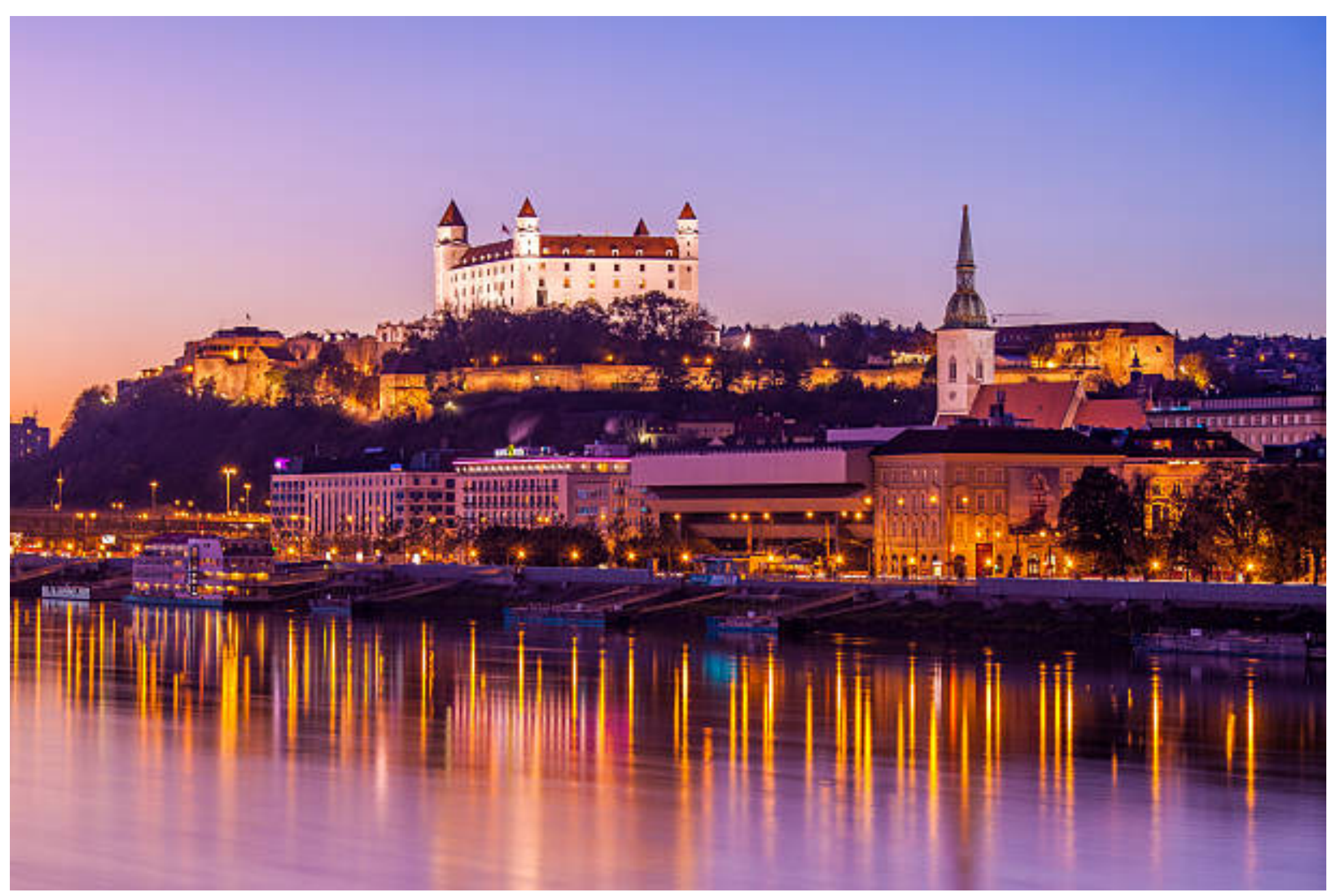

**Figure 1.** Bratislava, Slovakia, at dusk. The streetlamps and illuminated riverside buildings produce tall, parallel-sided pillars of light reflected in the uniformly rippled waters of the Danube

The answer lies in the interplay between two effects: (1) Rippled water surface has local surface slopes that can be described by a height function $z = \zeta(x, y)$ (the symbol is chosen to align with previous related analysis on sea waves[1]). Each small patch of the surface has a local tangent plane that acts as a tiny mirror. Only those patches whose slopes fall within a certain range can reflect the light source into the observer's eye since the specular reflection rule (angle of reflection of a ray is equal to the angle of incidence). (2) The eye or camera (represented optically as a pinhole, although it is a biconvex lens) projects the three-dimensional scene onto a two-dimensional (see clarification of this below) image plane. Points farther across the water are compressed by a factor proportional to $1/x$ in this process. This geometric compression counteracts the physical widening of the region of contributing facets.

The result is a reflection whose image width remains constant. However, we show that there is some deviation of this simple picture at both ends of the image.

Once reflection is introduced to physical-optical systems interesting images arise, as occurs with variously shaped mirrors[2-7]; so, the present analysis is part of the general theme of reflections in geometry and hence physical optics.

## 2. Previous incomplete analysis

Figure 2 shows a simplified side-on (elevation) view of the scene: the streetlight, the rippled water surface, and an observer. Earlier accounts of light-pillar formation correctly noted that a rippled water surface contains many small facets (tangent planes) acting as tiny mirrors[8]. A facet reflects the light source into the observer's eye only when its orientation satisfies the specular reflection rule. This explains why isolated bright patches appear on the water, but it does not explain why the *image* of these patches forms a tall, narrow, parallel-sided pillar.

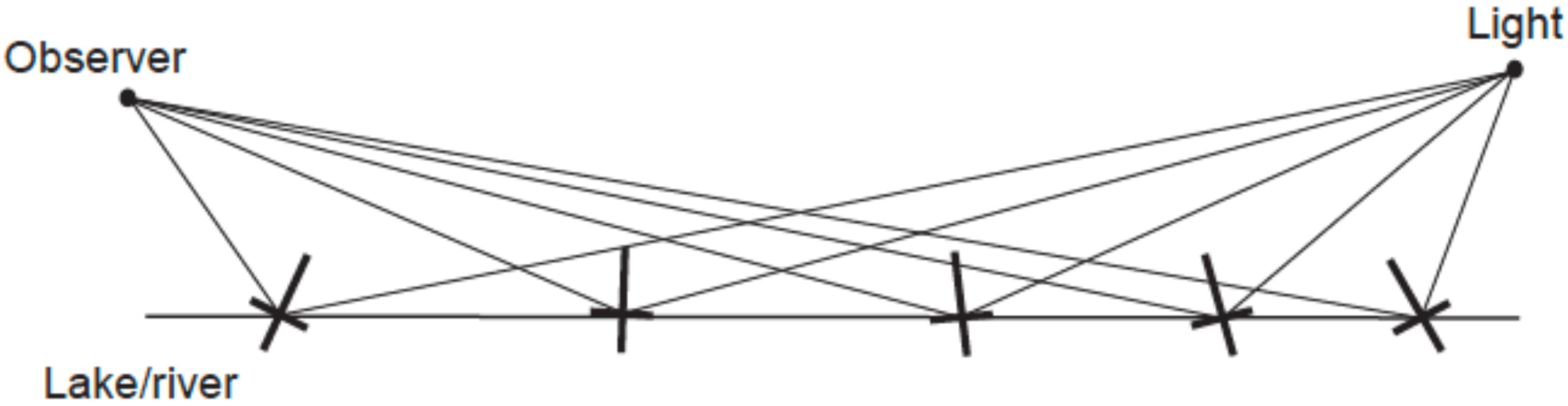

**Figure 2.** Side-on (elevation) schematic of the reflection geometry of a streetlight, the rippled water surface, and an observer[8]. Each bright patch on the water corresponds to a wave facet whose orientation satisfies the specular reflection rule (angle of reflection equals angle of incidence). The longer arm of each cross indicates the surface normal; the shorter arm shows the facet edge-on. This schematic explains why a single reflection *line* is visible, not simply a single spot like from a pure flat mirror; but it does not yet explain why the extended reflection forms a parallel-sided pillar in the image.

We show below that rays that reach the pinhole form an inverted image on that plane, and the geometry of this projection stretches the reflection vertically while preserving the parallelism of its sides. In other words, the pillar is the image-plane footprint of many specular patches arranged in the line of sight.

*Aside:* The optical projection through the eye's lens produces an inverted image on the retina that is not a flat plane but a section of a sphere, reversing both the vertical and horizontal axes (Fig. 3). This inversion is a consequence of geometric optics, not a neural operation. The brain does not rotate the image; instead, it rapidly interprets and re-maps the retinal signals – drawing on visuomotor and multisensory cues – so that our *perceived* world is upright and aligned with touch and action[9].

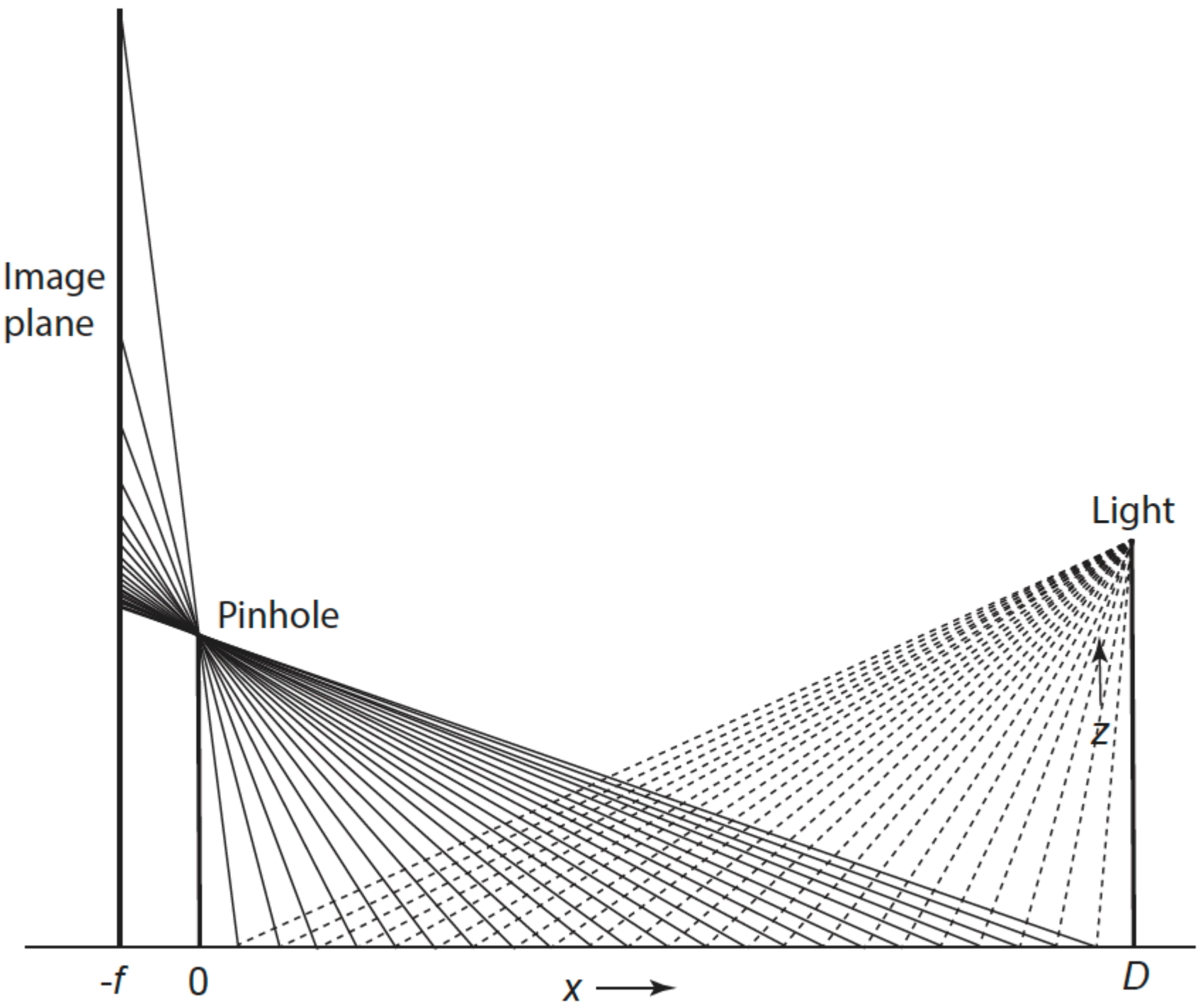

**Figure 3.** Elevation view of the pinhole projection system in which the $x$-axis represents the mean water surface in cross-section with the pinhole at $x = 0$ and the image plane at $x = -f$. Note the compression of light rays from larger $x$ values at the lower reaches of the *image*.

Once the pinhole model is written in coordinates, the appearance of long parallel-sided pillars follows directly from how changes in elevation angle (vertical tilts of the incoming rays) map to the vertical coordinate $v$ in the image plane. The remainder of this paper develops this interaction quantitatively.

## 3. Geometry of the Reflection

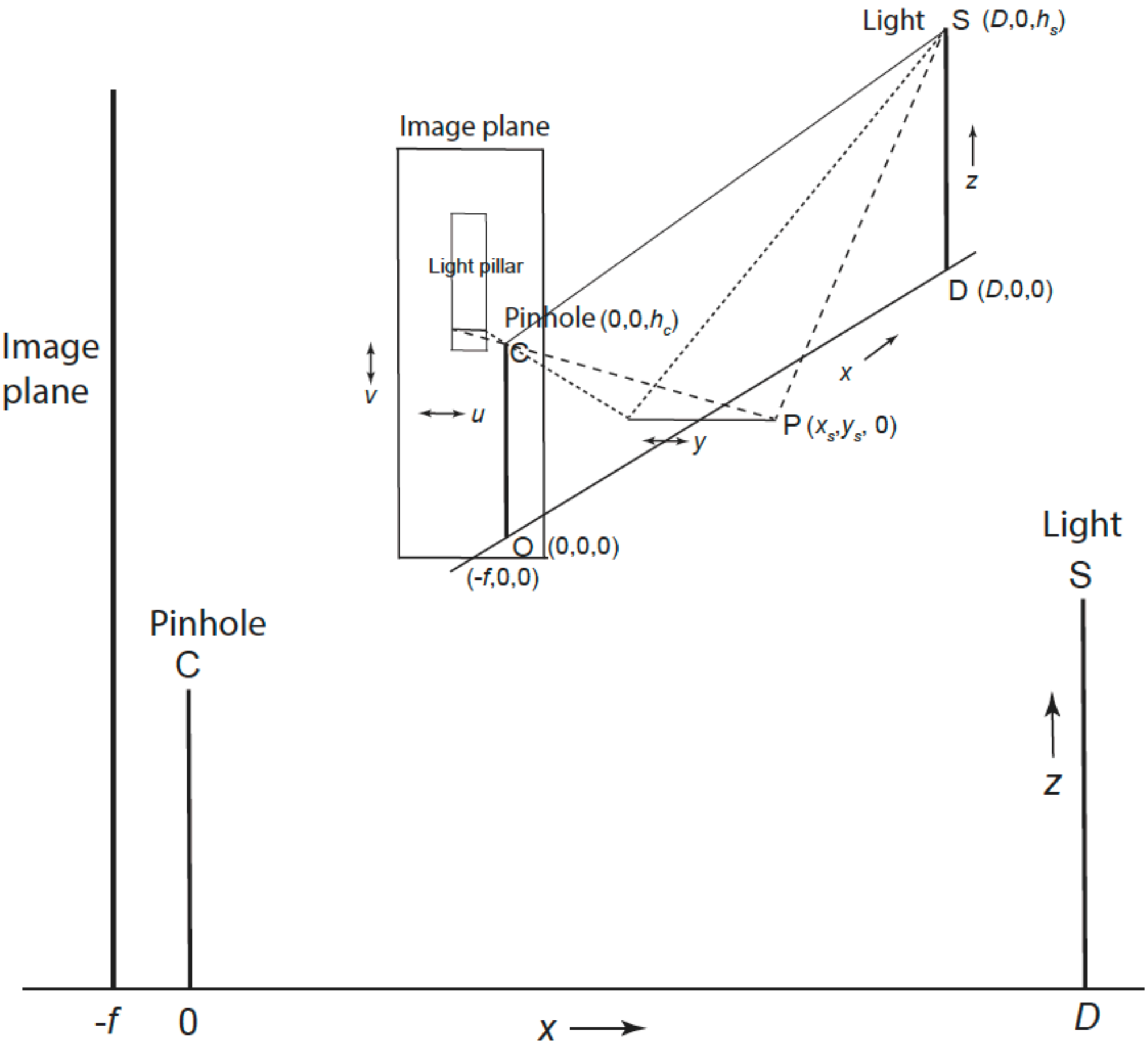

**Figure 4.** Geometry of the specular reflection. The camera pinhole is at $\mathrm{C} = (0, 0, h_c)$, the light source at $\mathrm{S} = (D, 0, h_\mathrm{s})$, and a surface point at $\mathrm{P}(x_\mathrm{s}, y_\mathrm{s}, 0)$. The rays to camera and source are $\mathbf{p} = (-x_\mathrm{s}, -y_\mathrm{s}, h_c)$ and $\mathbf{q} = (D - x_\mathrm{s}, -y_\mathrm{s}, h_\mathrm{s})$. Their unit directions are $\hat{\mathbf{u}} = \mathbf{u}/\| \mathbf{u} \|$and $\hat{\mathbf{v}} = \mathbf{v}/\| \mathbf{v} \|$. The bisector is $\mathbf{b} = \hat{\mathbf{u}} + \hat{\mathbf{v}}$, and the surface normal for a height field $z = \zeta(x, y)$ is $\mathbf{n} = (-\zeta_x, -\zeta_y, 1)$. The specular condition requires $\mathbf{n} \parallel \mathbf{b}$. Projection onto the plane $x = -f$ gives the image coordinates $(u, v)$. All vectors are shown in roman bold, with components and scalars in italics.

### 3.1 Coordinate System and Key Coordinates

We adopt a right-handed Cartesian coordinate system in which the undisturbed water surface lies in the plane $z = 0$. Horizontal coordinates are $(x, y)$, and height is $z$.
The three essential points in the geometry are:

$$\mathrm{C} = (0,0, h_\mathrm{c}),\ \mathrm{S} = (D, 0, h_\mathrm{s}),\ \mathrm{P} = (x, y, 0). \tag{1}$$

Here, C is the (camera) pinhole at height $h_\mathrm{c}$, S is the light source at height $h_\mathrm{s}$, and P is an arbitrary point on the water surface.

The vector from P to the pinhole is

$$\mathbf{p} = \mathrm{C} - \mathrm{P} = (-x,\ -y, h_{\mathrm{c}}), \tag{2}$$

and the vector from P to the source is

$$\mathbf{q} = \mathrm{S} - \mathrm{P} = (D - x,\ -y, h_{\mathrm{s}}). \tag{3}$$

Let the corresponding unit vectors be

$$\hat{\mathbf{p}} = \frac{\mathbf{p}}{\|\mathbf{p}\|}, \quad \hat{\mathbf{q}} = \frac{\mathbf{q}}{\|\mathbf{q}\|}. \tag{4}$$

### 3.2 Specular Reflection Condition

For specular reflection, the surface normal at P must bisect the angle between the incident and reflected rays. Thus the (un-normalised) bisector vector is

$$\mathbf{b} = \hat{\mathbf{p}} + \hat{\mathbf{q}}. \tag{5}$$

Appendix B proves that **b** bisects the two ray directions.

The specular condition is therefore that **n**, the surface normal at P is parallel (∥) to **b**:

$$\mathbf{n} \parallel \mathbf{b}. \tag{6}$$

### 3.3 Surface Slopes and the Normal Vector

A rippled water surface can be described locally by a *height field* as originally described by Cox and Munk[1]; and we use their $\zeta$ for this:

$$z = \zeta(x, y). \tag{7}$$

The tangent vectors of this height field are the partial derivatives of $\zeta$ with respect to *x* and *y* respectively[10]:

$$\mathbf{t}_x = (1, 0,\ \zeta_x),\ \ \mathbf{t}_y = (0, 1, \zeta_y). \tag{8}$$

So, the (un-normalised) surface normal is the cross-product of these two vectors:

$$\mathbf{n} = \mathbf{t}_x \times \mathbf{t}_y = (-\zeta_x, -\zeta_y, 1). \tag{9}$$

This result is derived in Appendix A.

If we define the slope components of a *surface facet* as $m_x$ and $m_y$ we see that the normal encodes these as well, because it contains these as the *x*- and *y*- components (albeit negative). Thus, the two projections on to the horizontal plane have the same magnitude, namely

$$m_x = \zeta_x, \ m_y = \zeta_y. \tag{10}$$

**3.4 Required Lateral Slope for Specular Reflection**

We now determine the slope $m_y$ that is required at a point $\mathrm{P} = (x, y)$ so that it will reflect the light from the source through the pinhole.

We restrict attention to the horizontal plane. The horizontal components of the ray directions are

$$\mathbf{q}_{hor} = (D - x, -y), \ \ \mathbf{p}_{hor} = (-x, -y). \tag{11}$$

The so called far-field regime that is relevant to light pillars applies as follows to the lateral displacement (where the vertical lines denote magnitude):

$$|\, y \,| \ll x, \qquad |\, y \,| \ll (D - x). \tag{12}$$

Normalising the ray directions yields the required lateral slope. Appendix E gives an alternative (perhaps more familiar to some readers) angular equivalent of this derivation. Thus,

$$m_y(x, y) \approx -\frac{y}{2}\left(\frac{1}{x} + \frac{1}{D-x}\right). \tag{13}$$

This is the key relationship: the required slope grows as the facet moves away from the central line $y = 0$, and it grows more rapidly near the far shore where $D - x$ becomes small but in the middle section of the body of water the dependence is effectively linear in $1/x$.

### 3.5 Lateral Acceptance Region

A real water surface cannot supply arbitrarily large slopes. There is a distribution of facet slopes that we will not consider in detail here (see Section 5); but it is plausible that there would be a maximum attainable slope magnitude that satisfies the specular reflection rule.

$$| m_y | \leq m_{\max}. \tag{14}$$

A point $\mathrm{P} = (x, y)$ contributes to the reflection only if the required slope (13) satisfies this bound:

$$\frac{y}{2}\left(\frac{1}{x} + \frac{1}{D-x}\right) \leq m_{\max}. \tag{15}$$

Solving for $| y |$ gives the physical half-width $W(x)$ of the contributing region:

$$W(x) = 2m_{\max} \frac{x(D-x)}{D}. \tag{16}$$

This width vanishes at $x = 0$ and $x = D$, and it is maximal at $x = D/2$, and grows quadratically with distance from the camera. Appendix E discusses this in more detail.

## 4. Projection and Image Formation

### 4.1 Pinhole Projection

The image plane (see Figs 3 and 4) is vertical at

$$x = -f, \tag{17}$$

with the pinhole at $\mathrm{C} = (0, 0, h_c)$.

The point $\mathrm{P} = (x, y, 0)$ projects to the image coordinates $(u, v)$ as is explained in detail in Appendix C:

$$u = -\frac{f\,y}{x+h_c}, \quad v = -\frac{f\,h_c}{x+h_c}. \tag{18}$$

For $x \gg h_c$, the horizontal coordinate simplifies to

$$u \approx -\frac{f\,y}{x}. \tag{19}$$

Thus, horizontal distances on the water shrink by a factor of $1/x$ in the image plane.
The physical half-width of the contributing region on the water surface is given by Eq. 16, so under the pinhole projection, this physical width is compressed by the factor $f/(x + h_c)$, giving the image-plane half-width of

$$u_{\max}(x) = \frac{f\,W(x)}{x+h_c}. \tag{20}$$

Substituting Eq. (16) yields

$$u_{\max}(x) = 2m_{\max}f\,\frac{x(D-x)}{D(x+h_c)}. \tag{21}$$

This expression is valid for the entire interval $0 < x < D$. It shows that the image width: (1) vanishes at the near-shore ($x = 0$); (2) vanishes at the far shore ($x = D$); (3) is broadest near the middle; and (4) narrows strongly as $x \to D,$ because the factor $D - x$ becomes small.
However, for most of the range where $x \gg h_c$, the denominator satisfies $x + h_c \approx x$, giving the simplified form

$$u_{\max}(x) \approx 2m_{\max}f\,\frac{D-x}{D}. \tag{22}$$

This approximation explains the parallel *appearance* of the pillar over much of its length, while Eq. (21) captures the pronounced taper close to the far shore. Meanwhile Eq. (22) shows that the image width varies linearly with $D - x$. Over most of the interval $0 < x < D$, this variation is weak, so the sides appear parallel and this is obviously a property of the projection, not of the water surface.

## 5. Brightness and Slope Statistics

### 5.1 Gaussian Distribution of Slopes

Surface slopes on wind-rippled water are well approximated by an isotropic Gaussian[1]:

$$p(m_x, m_y) = \frac{1}{2\pi\sigma^2} \exp(-\frac{m_x^2+m_y^2}{2\sigma^2}). \quad (23)$$

Thus, moderate slopes are common; steep slopes are exponentially rare.

### 5.2 Sparkled Edges

Near the geometric boundary of the light pillar, the required slope approaches the maximum attainable value $m_{\max}$. In this region the Gaussian distribution lies in its tails, so suitable facets are rare. Thus, the edges of the pillar appear as isolated bright points that are the familiar sparkled boundary. Each sparkle corresponds to a transient alignment of a rare steep facet.

### 5.3 Brightness Variation Along the Pillar

Photographs (*e.g.*, Figure 1) show that the image of the pillar is (perhaps) brightest near its base (reflections from near the far shore) and fades toward the top. This follows from the mapping in Eq. (18), $v(x) = fh_C/x$ because differentiating this with respect to $x$ gives the slope in the image intensity as a function of the density of the reflectors along the horizontal $x$-coordinate:

$$\frac{dv}{dx} = -\frac{fh_C}{x^2}. \quad (24)$$

Thus, a fixed vertical interval in the image corresponds to a physical interval proportional to $x^2$. In other words, as can be seen in Figure 3, farther water contributes more facets per unit image height than the closer water, producing greater brightness near the bottom of the raw image (near the top of the perceived pillar after pinhole inversion).

### 5.4 A competing Effect on Brightness

As just shown, projection compression comes from Eq. (18), and the facets near the far shore collect light from a larger physical area. This tends to increase brightness toward the horizon. But we must also consider the rarity of large slopes as expressed in Eq. (23) and the required slope increases as $x \to D$. In this limit large slopes are exponentially unlikely under Eq. (22). This tends to decrease brightness in the image towards the horizon. Overall, under typical surface conditions, statistical suppression dominates, so the pillar is brightest near the pinhole and fades toward the far shore; we posit that this effect is evident in Figure 1 but detailed image analysis would be needed to back-up this claim!

## 6. Conclusions

A rippled water surface acts as a mosaic of tiny mirrors. Specular reflection selects only those facets whose normals bisect the rays from source to pinole. The finite range of attainable slopes defines a physical contributing region whose width grows linearly with distance.

Pinhole projection compresses distant regions by $1/x$, cancelling (except at the ends) this geometric widening. The result is a tall, narrow, almost parallel-sided pillar, which is a property of imaging geometry, not of the water itself.

Brightness variations arise from the Gaussian statistics of surface slopes and the competing effects of projection compression and slope rarity. The same framework extends naturally to reflections of the Sun viewed from elevated positions.

## 7. Use of AI

Microsoft Copilot was used to refine and abbreviate the presentation, but not as final arbiter.

## 8. Acknowledgements

Dr Daniel Daners is thanked for valuable discussions.

## Appendix A - Surface Normal for a Height Field

A water surface can be described locally by a height function[1]

$$z = \zeta(x, y). \tag{A1}$$

A point on the surface is:

$$\mathrm{P}(x, y) = (x, y,\ \zeta(x, y)). \tag{A2}$$

The tangent vectors are obtained by differentiating P[10]:

$$\mathbf{t}_x = \frac{\partial P}{\partial x} = (1, 0, \zeta_x),\ \ \mathbf{t}_y = \frac{\partial P}{\partial y} = (0, 1, \zeta_y). \tag{A3}$$

The (un-normalised) surface normal is given by the cross product:

$$\mathbf{n} = \mathbf{t}_x \times \mathbf{t}_y = \begin{vmatrix} \hat{\mathbf{i}} & \hat{\mathbf{j}} & \hat{\mathbf{k}} \\ 1 & 0 & \zeta_x \\ 0 & 1 & \zeta_y \end{vmatrix} = in\ component\ form\ \left(-\zeta_x, -\zeta_y, 1\right). \tag{A4}$$

Hence, this is Eq. (9).

**Appendix B - Why the Bisector Vector Enforces Specular Reflection**

Let $\hat{\mathbf{p}}$ and $\hat{\mathbf{q}}$ be the unit vectors along the reflected and incident rays at a surface point. Define the bisector vector as

$$\mathbf{b} = \hat{\mathbf{p}} + \hat{\mathbf{q}}. \tag{B1}$$

Compute the dot products:

$$\mathbf{b} \cdot \hat{\mathbf{p}} = 1 + \hat{\mathbf{p}} \cdot \hat{\mathbf{q}}, \mathbf{b} \cdot \hat{\mathbf{q}} = 1 + \hat{\mathbf{p}} \cdot \hat{\mathbf{q}}. \tag{B2}$$

Since these are equal, the angles between **b** and each ray-direction are equal. Thus **b** bisects the angle between the incident and reflected rays, and the specular condition is Eq. (6):

$$\mathbf{n} \parallel \mathbf{b}. \tag{B3}$$

**Appendix C - Pinhole Projection Formulae**

The pinhole is located at

$$\mathrm{C} = (0, 0, h_{\mathrm{c}}), \tag{C1}$$

and the image plane is the vertical plane at

$$x = -f. \tag{C2}$$

A surface point is:

$$\mathrm{P} = (x, y, 0). \tag{C3}$$

The ray from C to P can be written as follows where the capital letters represent the Cartesian coordinates of the points:

$$\mathrm{C} + \lambda(\mathrm{P} - \mathrm{C}). \tag{C4}$$

The intersection with the plane at $x = -f$ occurs when

$$0 + \lambda(x - 0) = -f \Rightarrow \lambda = -\frac{f}{x}. \tag{C5}$$

Thus, the projected coordinates are:

$$u = \lambda(0 - y) = -\frac{f\,y}{x}, \tag{C6}$$

$$v = \lambda(0 - h_c) = -\frac{f\,h_c}{x}. \tag{C7}$$

These are the projection relations given in Eq. (18) and used beyond.

**Appendix D - Incidence-Plane Normal**

The incidence plane contains the three points C, P, and S whose coordinates can be interpreted as position vectors. Two non-parallel vectors in this plane are, therefore

$$\mathbf{u} = \mathrm{C} - \mathrm{P}, \quad \mathbf{v} = \mathrm{S} - \mathrm{P}. \tag{D1}$$

A normal to the incidence plane spanned by **u** and **v** is given by their cross-product[10]

$$\mathbf{n}_{\mathrm{inc}} = \mathbf{u} \times \mathbf{v}. \tag{D2}$$

This vector is useful for visualising the geometry of the reflection, but it is not required in the main derivation.

**Appendix E - Angular (Azimuthal) Formulation and Its Equivalence**

For completeness and synergy with an angular approach in optics, we show that the coordinate-based derivation of the required slope is equivalent to our derivation.
Let the horizontal coordinates of a surface point be written in polar form:

$$(x, y) = (p\cos\ \phi, p\sin\ \phi). \quad \text{(E1)}$$

The source and camera are at

$$\mathrm{S} = (D, 0, h_s),\ \ \mathrm{C} = (0, 0, h_c). \quad \text{(E2)}$$

In the far-field regime, the horizontal projections of the ray directions are:

$$(D - x,\ -y),\ (-x,\ -y). \quad \text{(E3)}$$

The corresponding azimuths are:

$$\phi_s \approx -\frac{y}{D-x},\quad \phi_c \approx -\frac{y}{x}. \quad \text{(E4)}$$

The bisector (bis) azimuth is:

$$\phi_{\text{bis}} = \frac{\phi_s + \phi_c}{2} \approx -\frac{y}{2}\left(\frac{1}{x} + \frac{1}{D-x}\right). \quad \text{(E5)}$$

For a height field $z = \zeta(x, y)$, the horizontal projection of the surface normal is readily seen to be proportional to $(\zeta_x, \zeta_y)$. Thus, the lateral slope required for specularity is:

$$m_y(x, y) \approx -\frac{y}{2}\left(\frac{1}{x} + \frac{1}{D-x}\right), \quad \text{(E6)}$$

which is identical to the coordinate-based result in Eq. 13, thus confirming the equivalence of the two approaches.

## 9. References

1. Cox CD, Munk W. Measurement of the roughness of the sea surface from photographs of the sun's glitter. *Journal of the Optical Society of America.* 1954;44(11):838-850.
2. Kuchel PW. Anamorphoscopes: a visual aid for circle inversion. *The Mathematical Gazette.* 1979;63(424):82-89.
3. Hickin P. Anamorphosis. *The Mathematical Gazette.* 1992;76(476):208-221.
4. Kuchel PW. Is there light at the end of the funnel? *The Mathematical Gazette.* 1994;78(483):336-339.
5. Sharp J, Nickel BG, Hunt JL. Anamorphoscopes: another look at circle-inverting mirrors. *The Mathematical Gazette.* 2011;95(532):1–16.
6. Berry MV. Inflection reflection: images in mirrors whose curvature changes sign. *European Journal of Physics.* 2021;42:1-12.
7. Berry MV. Distorted mirror images organised by cuspoid and umbilic caustics. *Journal of Optics.* 2021;23:125402.
8. Minnaert MGJ. *Light and Color in the Outdoors.* New York: Srpinger-Verlag; 1974.
9. Kandel ER, Schwartz JH, Jessell TM, Siegelbaum SA, Hudspeth AJ. *Principles of Neural Science.* 5th ed. New York: McGraw-Hill; 2013.
10. Gray A, Abbena E, Salamon S. *Modern Differential Geometry of Curves and Surfaces with Mathematica.* 3rd ed. Boca Raton: Chapman & Hall/CRC; 2006.